\documentclass[12pt]{amsart}

\usepackage{amsmath,amsthm,amsfonts,amscd,flafter,epsf}

\newcommand\commentable[1]{#1}

\newcounter{bean}

\newtheorem{theorem}{Theorem}[section]
\newtheorem{prop}[theorem]{Proposition}

\newtheorem{remark}[theorem]{Remark}

\def\endproof{\relax\ifmmode\expandafter\endproofmath\else
  \unskip\nobreak\hfil\penalty50\hskip.75em\hbox{}\nobreak\hfil\bull
  {\parfillskip=0pt \finalhyphendemerits=0 \bigbreak}\fi}
\def\endproofmath$${\eqno\bull$$\bigbreak}
\def\bull{\vbox{\hrule\hbox{\vrule\kern3pt\vbox{\kern6pt}\kern3pt\vrule}\hrule}}

\newcommand{\R}{\mathbb{R}}

\newcommand{\Z}{\mathbb{Z}}

\newcommand{\CP}[1]{{\mathbb{CP}}^{#1}}
\newcommand{\CPbar}{{\overline{\mathbb{CP}}}^2}
\newcommand{\Zmod}[1]{\Z/{#1}\Z}

\newcommand{\ModSW}{\ModSWfour}

\newcommand{\ModSWthree}{\mathcal{N}}

\newcommand{\ModSWfour}{\mathcal{M}}
\newcommand{\ModFlow}{\ModSWfour}
\newcommand{\UnparModFlow}{\widehat\ModFlow}
\newcommand{\SW}{SW}

\newcommand{\SWConfig}{\mbox{${\mathcal B}$}}
\newcommand{\SWConfigIrr}{\mbox{$\SWConfig^*$}}

\newcommand{\SpinC}{\Spin^{c}}
\newcommand{\spinct}{\mathfrak{t}}
\newcommand{\Spin}{{\mathrm{Spin}}}

\newcommand{\goesto}{\mapsto}

\newcommand{\CSD}{\mbox{${\rm CSD}$}}

\newcommand\Wedge{\Lambda}

\newcommand\abuts\Rightarrow
\newcommand\Sym{\mathrm{Sym}}

\newcommand\Alg{\mathbb{A}}

\newcommand\Jac{\mathcal J}
\newcommand\CritMan{C}

\newcommand\Cinfty{C^{\infty}}

\newcommand{\PD}{\mathrm{PD}}

\newcommand{\ModSWThree}{\ModSWthree}

\newcommand\spinc{\mathfrak s}
\newcommand\XExt{X^{+}}

\newcommand\mult{e}
\newcommand\glue{\gamma}
\newcommand\Restrict[1]{\rho_{_{#1}}}
\newcommand\ProdSpinCs{\spinc'\in\spinc_{1}\#\spinc_{2}}

\title{On Embedding Circle-Bundles in Four-Manifolds}
\author[Peter Ozsv{\'a}th]{Peter Ozsv\'ath} 
\thanks{The first author was partially supported by NSF grant number 
DMS 9971950.}
\address{Department of 
Mathematics, Michigan State University,
East Lansing, Michigan 48824}

\author[Zolt{\'a}n Szab{\'o}]{Zolt{\'a}n Szab{\'o}}
\thanks{The second author was partially 
supported by NSF grant number DMS 970435 and a
Sloan Fellowship} 
\address{Department of Mathematics, University of Michigan,
Ann Arbor, Michigan 48109-1109}

\begin{document}

\begin{abstract} 
    In this paper, we demonstrate an obstruction to 
    finding certain
    splittings of 
    four-manifolds along sufficiently twisted circle bundles over Riemann 
    surfaces, arising from Seiberg-Witten theory. These obstructions 
    are used to show a non-splitting result for algebraic 
    surfaces of general type.
\end{abstract}

\maketitle
\section{Introduction}

\label{sec:Introduction}
Let $Y(n,g)$ denote the circle bundle over a genus $g$ surface with
Euler number $n$.
Our main result in this paper is the following:

\begin{theorem}
\label{thm:GenType}
If $X$ is a complex surface of general type, and $Y(n,g)$ is
circle-bundle over a Riemann surface of genus $g$, whose Euler number
$n$ satisfies $|n|\geq 2g-1$, then $X$ admits no splittings along an
embedded copy of $Y=Y(n,g)$ of the form $X=X_1\#_Y X_2$ with
$b_2^+(X_1), b_2^+(X_2)>0$.
\end{theorem}

In the above theorem, the quantity $b_{2}^{+}(Z)$ of a four-manifold
$Z$ with boundary denotes the maximal dimension of a positive-definite
subspace for the intersection form on $H^2(Z,\partial Z;\R)$. It is
suggestive to compare the hypothesis that $|n|\geq 2g-1$ with the
``adjunction inequality'' for surfaces of non-negative square
(see~\cite{KMthom} or
\cite{MSzT}). Indeed, the hypotheses of Theorem~\ref{thm:GenType} are
sharp: if we allow $b_2^+(X_2)=0$ or $|n|\leq 2g-2$, there are many
examples of such splittings, obtained by blowing up smoothly embedded
complex curves $C$ in $X$, and splitting $X$ along the boundary of the
tubular neighborhood of $C$.

Moreover, the situation for elliptic surfaces is quite different, 
as we see below:

\begin{theorem}
    \label{thm:ESurf}
     \begin{list}
{(\arabic{bean})}{\usecounter{bean}\setlength{\rightmargin}{\leftmargin}}
\item
\label{item:EulerOne}
Every simply-connected elliptic surface with $b_2^+(X)>3$ admits
a splitting along $Y(1,1)$ with $b_2^+(X_i)>0$. 
\item
\label{item:EulerGen}
For each $n>0$, there is a simply-connected elliptic surface $X$
which admits a splitting along $Y(n,1)$ with $b_2^+(X_i)>0$.
\end{list}
\end{theorem}

\begin{proof}
    A splitting of Type~(\ref{item:EulerOne}) is given as follows. 
    Note
    that
    $Y(1,1)$ is the mapping cylinder of a (single) Dehn twist on the
    torus. Thus, if we begin with the rational elliptic surface
    $$\pi\colon E(1)=\CP{2}\# 9 \CPbar\longrightarrow\CP{1},$$ let
    $x\in \CP{1}$ denote the image of a fishtail fiber, and let
    $\Delta$ be a disk around $x$ containing no other singular points
    for $\pi$, then $\pi^{-1}(\partial \Delta)=Y(1,1)$ splits $E(1)$
    into a pair of elliptic fibrations $Z_{1}$ and $Z_{2}$ over disks.
    Thus we can realize $E(3)$ as a union of fiber sums $E(1)\# Z_{1}$
    and $Z_{2}\# E(1)$ joined along $Y(1,1)$, where $E(3)$ is the
    fiber sum of three copies of $E(1)$. Neither side is 
    negative-definite: both sides contain a torus of 
    square zero and a sphere (constructed from vanishing cycles) 
    which meets this torus in a single, 
    positive point. 
    Since every simply-connected
    elliptic surface with $b_2^+\geq 3$ can be obtained from $E(3)$ by
    fiber sums with $E(1)$, logarithmic transformations, and blow-ups
    (see~\cite{Moishezon}, or~\cite{FriedMorgBook})
    the result follows.
    
    A splitting of Type~(\ref{item:EulerGen}) is realized by finding
    an elliptic surface $Z$ over $\CP{1}$ which contains $n$ singular
    values for the elliptic fibration whose holonomy is a Dehn
    twist along a given curve in the fiber. In fact, it is a theorem
    of Moishezon (see~\cite{Moishezon}, also Theorem~3.6 in Chapter~2
    of~\cite{FriedMorgBook}) that if $Z$ is a nodal elliptic surface
    without multiple fibers and $2m$ singular fibers, then we can
    think of the monodromy representation around $m$ of the singular
    fibers, of which we select $n$, as being a ($+1$) Dehn twist
    around a fixed non-separating curve in the fiber, and the
    monodromy around the remaining $m$ as being a Dehn twist around
    another curve.  Let $\Delta$ be a disk in $\CP{1}$ which contains
    only the $n$ distinguished singular points and no others.  
    Now, it is easy to see that $\pi^{-1}(\partial
    \Delta)=Y(n,1)$, which separates the elliptic surface. Forming
    fiber sums with rational elliptic surfaces on both sides as
    before, we get a decomposition of the elliptic surface $E(1)\# Z\#
    E(1)$ along $Y(n,1)$ into two pieces with $b_2^+>0$.
    \end{proof}
    
\begin{remark}
    Note that the hypothesis that $b_{2}^{+}(Z)>3$ above is necessary:
    the elliptic surface $E(2)$ admits
    no decomposition along $Y(n,g)$ with $g\geq 2n-1$ and 
    $b_{2}^{+}(X_{i})>0$.
    This follows from the fact that $E(2)$ has a single basic class, 
    together with the vanishing result, 
    Theorem~\ref{thm:VanishingTheorem}, from Section~\ref{sec:Proof}.
\end{remark}
        
\begin{remark}
Using the above decomposition (Type~\ref{item:EulerGen}) as a building
block, it is possible to construct symplectic four-manifolds $Z$ which
decompose along $Y(n,g)$ with $n$ and $g$ arbitrarily large, such that
both sides have positive $b_2^+$. For example, one can start with an
elliptic surface $X$ decomposed along $Y(n,1)$ in the manner of
Theorem~\ref{thm:ESurf}, and find a symplectic torus $T\subset X$
(which is symplectic for a form arbitrarily close to a K\"ahler form
for $X$) which meets $Y(n,1)$ in a fiber circle for the Seifert
fibration of $Y(n,1)$, and has square zero.  Forming the fiber sum of
$X$ with, say, $T^2\times\Sigma_{g-1}$ (by gluing $T\subset X$ to
$T^2\times p$), we obtain $Z$ as claimed.
\end{remark}
    
Theorem~\ref{thm:GenType} follows from a ``vanishing theorem,''
Theorem~\ref{thm:VanishingTheorem}, according to which a certain sum
of Seiberg-Witten invariants for $X$ vanishes whenever $X$ splits into
two pieces with $b_2^+(X_i)>0$ along $Y(n,g)$, when $|n|\geq
2g-1$. This is a more refined vanishing statement than the usual
vanishing theorem over $S^3$: in particular there are manifolds with
non-trivial Seiberg-Witten invariants which admit such splittings, as
is illustrated by Theorem~\ref{thm:ESurf}. The vanishing theorem is
proved by looking at the ends of the moduli spaces of flows to the
reducibles: this is also the philosophy adopted by Austin and Braam
in~\cite{AusBra}, see also~\cite{StipSzabo}. In the case where $g=1$,
it is interesting to compare the vanishing theorem with a certain
vanishing theorem for Donaldson polynomials proved by Morgan, Mrowka,
and Ruberman (Theorem~16.0.1 of~\cite{MMR}).

We will give the proof of Theorem~\ref{thm:GenType} in 
Section~\ref{sec:Proof}, after stating and proving the more general 
vanishing result on which it is based.

\section{The Vanishing Theorem}
\label{sec:Proof}

To state the vanishing theorem which implies
Theorem~\ref{thm:GenType}, we must introduce some notation. We think
of the Seiberg-Witten invariant of a smooth, oriented, closed
four-manifold $X$ (with a ``homology orientation'' -- an orientation
on $(H^{0}\oplus H^{1}\oplus H^{+})(X;\R)$) and $\SpinC$ structure
$\spinc$ as a homogenous polynomial map $$\SW_{(X,\spinc)}\colon
\Alg(X)\longrightarrow \Z$$ of degree
$$d(\spinc)=\frac{c_{1}(\spinc)^{2}-(2\chi+3\sigma)}{4}$$ on the
algebra $$\Alg(X)=\Z[U]\otimes_{\Z}\Wedge^{*}H_{1}(X;\Z),$$ where $U$
is a two-dimensional generator, and $\Wedge^* H_1(X;\Z)$ is the
exterior algebra on the first homology of $X$ (graded in the obvious
manner). This algebra maps surjectively to the cohomology ring of the
irreducible configuration space $\SWConfigIrr(X,\spinc)$ of pairs
$[A,\Phi]$ of $\SpinC$ connections $A$ and somewhere non-vanishing
spinors $\Phi$ modulo gauge. (We denote the full configuration space
of pairs modulo gauge, i.e. where $\Phi$ is allowed to vanish, by
$\SWConfig(X,\spinc)$.)  As usual, the Seiberg-Witten invariant is
obtained by cohomological pairings of these cohomology classes with the
fundamental cycle of the moduli space $\ModSW(X,\spinc)$ of solutions
to the Seiberg-Witten equations, which is naturally induced from the
homology orientation.

As in Section~\ref{sec:Introduction}, let $Y=Y(n,g)$ be the circle bundle over 
a Riemann surface with Euler number $n$ over a surface $\Sigma$ of 
genus $g>0$. Throughout this section, we assume that
$$|n|\geq 2g-1.$$ 

Recall that $H^{2}(Y;\Z)\cong \Z^{2g}\oplus (\Zmod{n})$, where 
the $\Zmod{n}$ factor is generated by multiples of the pull-back $\pi^{*}$
of the 
orientation class of $\Sigma$. Indeed, there is a canonical $\SpinC$ 
structure $\spinct_{0}$ over $Y$ associated to the two-plane field orthogonal to the 
circle directions. Thus, forming the tensor product with $\spinct_{0}$
gives a canonical identification
$$\SpinC(Y)\cong H^{2}(Y;\Z).$$ In particular, there are $n$ 
distinguished $\SpinC$ structures $\spinct_{e}$ over $Y$, indexed by 
$e\in \Zmod{n}$ (thought of as a subset 
of $H^{2}(Y;\Z)$).

\begin{theorem}
    \label{thm:VanishingTheorem}
    Let $X$ be a smooth, closed, oriented four-manifold which splits 
    along an embedded copy of $Y=Y(n,g)$ with $|n|>2g-1$, so that 
    $X=X_{1}\#_{Y}X_{2}$ with $b_{2}^{+}(X_{i})>0$ for $i=1,2$. 
    Fix a $\SpinC$ structure $\spinc$ on $X$, and let
    $\spinc|_{Y}=\spinct$. If $\spinct$ is not one of the $n$ 
    distinguished $\SpinC$ structures on $Y$, then 
    $\SW_{(X,\spinc)}\equiv 0$. Similarly, if $\spinct=\spinct_{\mult}$
    for $2g-2<\mult<n$, then
    $\SW_{(X,\spinc)}\equiv 0$.
    Otherwise, if $\spinct=\spinct_{\mult}$ for $i=0,...,2g-2$, we have that
    \begin{equation}
	\label{eq:VanishingTheorem}
	\sum_{\{\spinc'|\spinc'|_{X_{1}}=\spinc|_{X_{1}}, 
    \spinc'|_{X_{2}}=\spinc|_{X_{2}}\}}\SW_{(X,\spinc')}\equiv 0.
    \end{equation}
\end{theorem}

Note that the 
inclusion $Y\subset X$ gives rise to a coboundary map $\delta\colon 
H^{1}(Y;\Z)\rightarrow H^{2}(Y;\Z)$, whose image we denote by $\delta 
H^{1}(Y;\Z)$. Another way of stating 
Equation~\eqref{eq:VanishingTheorem} is:
$$\sum_{\eta\in \delta H^{1}(Y;\Z)}\SW_{(X,\spinc+\eta)}\equiv 0.$$

The above theorem is proved by considering the ends of certain moduli 
spaces over cylindrical-end manifolds. In general, these ends are 
described in terms of the moduli spaces of the boundary $Y$, and the 
moduli spaces of solutions on the cylinder $\R\times Y$ (using a 
product metric and perturbation).

Specifically, let $Y$ be a three-manifold, and let $\ModSWThree_{Y}(\spinct)$ denote 
the moduli space of solutions to the three-dimensional Seiberg-Witten 
equations over $Y$ in the $\SpinC$ structure $\spinct$. Given a pair 
of components $C_{1}$, $C_{2}$ in $\ModSWThree_{Y}(\spinct)$, let 
$\ModSW(C_{1},C_{2})$ denote the moduli space of solutions $[A,\Phi]$ 
to the Seiberg-Witten equations on $\R\times Y$ for which
\begin{eqnarray*}
    \lim_{t\goesto -\infty}[A,\Phi]|_{\{t\}\times Y}\in C_{1}
    &{\text{and}}&
    \lim_{t\goesto \infty}[A,\Phi]|_{\{t\}\times Y}\in C_{2}
\end{eqnarray*}
The theory of~\cite{MMR} can be adapted to give the moduli space 
$\ModSW(C_{1},C_{2})$ a Fredholm deformation theory, and a pair of 
continuous ``boundary value maps'' for $i=1,2$
$$\rho_{_{C_{i}}}\colon \ModFlow(C_{1},C_{2})\longrightarrow C_{i}$$
characterized by 
\begin{eqnarray*}
    \rho_{_{C_{1}}}[A,\Phi]=\lim_{t\goesto -\infty}[A,\Phi]|_{\{t\}\times Y}
    &{\text{and}}&
    \rho_{_{C_{2}}}[A,\Phi]=\lim_{t\goesto +\infty}[A,\Phi]|_{\{t\}\times Y}.
    \end{eqnarray*}
The moduli space $\ModFlow(C_{1},C_{2})$
admits a translation action by $\R$, and we let 
$\UnparModFlow(C_{1},C_{2})$ denote the quotient of this space by this action. 
The boundary value maps are $\R$-invariant, and hence 
induce boundary value maps on the quotient
$$\rho_{_{C_{i}}},\colon \UnparModFlow(C_{1},C_{2})\longrightarrow C_{i}.$$

As in~\cite{KMthom} (by analogy with the cases considered by Floer, 
see for instance~\cite{Floer}),
the solutions to the three-dimensional Seiberg-Witten equations are 
the critical points for a ``Chern-Simons-Dirac'' functional $\CSD$ 
defined on the configuration space $\SWConfig(Y,\spinct)$. The 
Seiberg-Witten equations on $\R\times Y$ can be naturally identified 
with upward gradient flowlines for this functional. (Strictly 
speaking, the functional $\CSD$ is real-valued only when the first 
Chern class $c_{1}(\spinct)$ is 
torsion; otherwise it is circle-valued.)

Solutions in $\ModSWThree(Y,\spinct)$ whose spinor vanishes identically 
correspond to flat connections on the determinant line bundle for 
$\spinct$. By 
analogy with the Donaldson-Floer theory, 
these solutions are usually called {\em reducibles}, and those with 
somewhere non-vanishing spinor are called {\em irreducibles}.

In the case where $Y$ is a non-trivial circle bundle over a Riemann
surface these moduli spaces were studied in~\cite{MOY}, see
also~\cite{HigherType}
(where $Y$ is endowed with a circle-invariant metric and the 
Seiberg-Witten equations over it 
are suitably perturbed).  Specifically, there is the following result:

\begin{theorem}
    \label{thm:ModY}
    Let $Y$ be a circle bundle over a Riemann surface with genus 
    $g>0$, and Euler number $|n|>2g-2$ (oriented as circle bundle 
    with negative Euler number). Then, the moduli space 
    $\ModSWThree_{Y}(\spinct)$ is empty unless $\spinct$ corresponds to 
    a torsion class in $H^{2}(Y;\Z)$. Suppose that 
    $\spinct=\spinct_{e}$ for $e\in\Zmod{n}\subset H^{2}(Y;\Z)$. Then
\begin{list}
{(\arabic{bean})}{\usecounter{bean}\setlength{\rightmargin}{\leftmargin}}
\item
\label{item:SmallK}
If $0\leq \mult \leq g-1$ then $\ModSWThree_Y(\spinct)$ contains two
components, a reducible one $\Jac$, identified with the Jacobian torus
$H^1(\Sigma;S^{1})$, and a smooth irreducible
component $\CritMan$ diffeomorphic to $\Sym^\mult(\Sigma)$. Both of these
components are non-degenerate in the sense of Morse-Bott. There is an
inequality $\CSD(\Jac)>\CSD(\CritMan)$, so
the space
$\UnparModFlow(\Jac,\CritMan)$ is empty.  The space
$\UnparModFlow(\CritMan,\Jac)$ is smooth of expected dimension $2\mult$;
indeed it is diffeomorphic to $\Sym^\mult(\Sigma)$ under the 
boundary value map
$$\rho_{_{\CritMan}}\colon\UnparModFlow(\CritMan,\Jac)\longrightarrow 
\CritMan\cong \Sym^{e}(\Sigma).$$
\newline
\item If $g-1< \mult\leq 2g-2$, the Seiberg-Witten moduli spaces over both
$Y$ and $\R\times Y$
in this $\SpinC$ structure are naturally identified with the	
corresponding moduli spaces in the $\SpinC$ structure $2g-2-\mult$, which
we just described. 
\newline
\item For all other $\mult$, $\ModSWThree_Y(\spinct)$ contains only
reducibles. Furthermore, it is smoothly identified with the Jacobian torus
$\Jac$.
\end{list}
\end{theorem}

When $\mult\neq g-1$, the above theorem is a special case of 
Theorems~1 and 2 of ~\cite{MOY} (see especially Corollary~1.5 
of~\cite{MOY}). When $\mult=g-1$, the case considered in that paper is 
not ``generic''. In fact, there is a natural perturbation (by some 
small multiple of the connection $1$-form of the Seifert fibration), 
which achieves the genericity stated above. This perturbation was used 
in~\cite{HigherType} to prove strong ``adjunction inequalities'' for 
manifolds which are not of simple type, and the above theorem in the 
case where $\mult=g-1$ 
is precisely Theorem~8.1 of~\cite{HigherType}. 
Note that the hypothesis $n>2g-2$ is required to separate the irreducible 
manifolds into distinct $\SpinC$ structures. Note also that if 
the orientation on $Y$ 
is reversed, the flow-lines reverse direction.

The proof of Theorem~\ref{thm:VanishingTheorem}
is obtained by considering the ends of the moduli 
spaces $\ModSW(X_{1},\spinc_{1},\Jac)$ of
Seiberg-Witten monopoles over the cylindrical-end manifold
$$\XExt_{1}=X_{1}\cup_{\partial X_{1}=\{0\}\times Y}[0,\infty)\times 
Y$$ 
in the $\SpinC$ structure $\spinc_{1}$, whose boundary values are 
reducible.
We will assume, as in that theorem, that 
$b_2^+(X_1)>0$. 
In general, moduli spaces of finite energy solutions to 
the Seiberg-Witten equations on a manifold with cylindrical ends are 
not compact. (The ``finite energy condition'' in this context
is equivalent to 
the hypothesis that the pair $[A,\Phi]$ has a well-defined boundary 
value.)
They do, however, have ``broken flowline'' 
compactifications (see~\cite{MMR} and \cite{Floer}). In particular, if $C$ is a 
component of $\ModSWThree(Y,\spinc_{1}|_{Y})$, then for generic 
perturbations, the moduli space
$\ModSW(X_{1},\spinc_{1},C)$ is a smooth manifold with finitely many 
ends
indexed by components $C_{1},\ldots,C_{n}$ in the moduli space
$\ModSWThree(Y,\spinc_{1}|_{Y})$, with 
$\CSD(C_{1})<\CSD(C_{2})<\ldots<\CSD(C_{n})<\CSD(C)$.
When all the $C_{i}$ are non-degenerate in the sense of Morse-Bott, 
and consist of irreducibles, a neighborhood of the corresponding end 
is diffeomorphic to the fibered product
$$\ModSW(X_{1},\spinc_{1},C_{1})\times_{C_1} \ModFlow(C_{1},C_{2})\times_{C_2}\ldots\times_{C_n}
\ModFlow(C_{n},C),$$
under a certain gluing map (provided that this space is a manifold -- 
i.e. provided that the various boundary value maps are transverse). 

In particular, suppose $X_{1}$ is an oriented four-manifold with 
boundary, whose boundary 
$\partial X_{1}$ is identified with $Y=Y(n,g)$ with the orientation 
described in Theorem~\ref{thm:ModY}. Then, it follows from that 
theorem that
if $\spinc_{1}|_{Y}=\spinct_{e}$ for $0\leq \mult\leq 2g-2$, then
$\ModSW(X_{1},\spinc_{1},\CritMan)$ is compact (since there are no ``intermediate'' 
critical manifolds to be added), and $\ModSW(X_{1},\spinc_{1},\Jac)$ 
has  a single end whose neighborhood is 
diffeomorphic to 
$$\ModSW(X_{1},\spinc_{1},\CritMan)\times (0,\infty).$$
(We use here the fact that  the restriction map 
$\rho_{_{\CritMan}}\colon\UnparModFlow(\CritMan,\Jac)\longrightarrow 
\CritMan$ is a diffeomorphism.) This gluing map 
$$\gamma\colon \ModSW(X_{1},\spinc_{1},\CritMan)\times (0,\infty)\longrightarrow 
\ModSW(X_{1},\Jac)$$ 
is compatible with 
restriction to compact subsets of $\XExt_{1}$; e.g. if we consider the 
compact subset $X_{1}\subset \XExt_{1}$, then
$$\lim_{T\goesto 
\infty}\gamma([A,\Phi],T)|_{X_{1}}=[A,\Phi]|_{X_{1}}.$$

We make use of the end of $\ModSW(X_{1},\spinc_{1},\Jac)$
in the following proposition.
Recall that the moduli space 
$\ModSW(X_{1},\spinc_{1},\CritMan)$ is a smooth, compact submanifold of 
the irreducible configuration space of $\XExt_{1}$.
It has a canonical 
top-dimensional homology class, denoted $[\ModSW(X_{1},\spinc_{1},\CritMan)]$,
induced from the ``homology orientation'' of $X_{1}$. It inherits 
cohomology classes by pulling back via the
boundary value map
$$\Restrict{\CritMan}\colon \ModSW(X_{1},\spinc_{1},\CritMan)$$
and from the 
natural map
$$i_{X_{1}}\colon \ModSW(X_{1},\spinc)\longrightarrow 
\SWConfigIrr(X_{1},\spinc)$$
given by restricting the pair $[A,\Phi]$ to the compact subset 
$X_{1}\subset\XExt_{1}$ (this restriction is irreducible from the 
unique continuation theorem for the Dirac operator).
The pairings with these classes can be 
thought of as a  ``relative Seiberg-Witten'' invariant. 

\begin{prop}
    \label{prop:TrivialRelativeInvariant}
    Suppose $b_2^+(X_1)>0$.
    Given any cohomology classes $a\in 
    H^{*}(\SWConfigIrr(X_{1},\spinc_{1}))$ 
    and $b\in H^{*}(\CritMan)$, the homology 
    class $[\ModSW(X_{1},\spinc_{1},\CritMan)]$ pairs 
    trivially with the class $i_{X_{1}}^{*}(a)\cup \rho_{_{C}}^{*}(b)$.
\end{prop}
    
\begin{proof}
    
    First, we reduce to the case where $b$ is absent (i.e. 
    zero-dimensional). This is done in two steps, first establishing 
    an inclusion
    \begin{equation}
	\label{eq:Inclusion}
	(i_{Y}\circ \Restrict{C})^{*}H^{*}(\SWConfigIrr(Y,\spinct))
    \subset i_{X_{1}}^{*}(\SWConfigIrr(X_{1},\spinc_{1})),
    \end{equation}
    where both are
    thought of as subsets of $H^{*}(\ModSW(X_{1},\spinc_{1},\CritMan))$,
    and then seeing that the map
    $$i_{Y}^*\colon H^{*}(\SWConfigIrr(Y,\spinct))\longrightarrow H^{*}(C)
	$$
    is surjective.
    
    To see Inclusion~\eqref{eq:Inclusion} we describe the geometric
    representatives for the generators of the cohomology ring
    $$H^{*}(\SWConfigIrr(Y,\spinct))\cong \Z[U]\otimes_{\Z}
    \Wedge^{*}(H_{1}(Y;\Z)).$$ Given a point $y\in Y$ and a line
    $\Lambda_{y}\subset W_{y}$ in the fiber of the spinor bundle over
    $y$, the class $U$ is Poincar\'e dual to the locus
    $V_{(y,\Lambda_{y})}$ of pairs $[B,\Psi]$ with
    $\Psi_{y}\in\Lambda_{y}$. Moreover, given a curve $\gamma\colon
    S^{1}\rightarrow X$, the corresponding one-dimensional cohomology
    class determined by the homotopy type of the map
    $$h_{\gamma}\colon \SWConfigIrr(Y,\spinct)\longrightarrow S^{1}$$
    given by measuring the holonomy of $B$ (relative to some fixed
    reference connection $B_{0}$) around $\gamma$; i.e. it is
    Poincar\'e dual to the preimage $V_{\gamma}$ of a regular value of
    $h_{\gamma}$.  This cohomology class is denoted $\mu[\gamma]\in
    H^{1}(\SWConfigIrr(Y,\spinct))$.  (Note that geometric
    representatives cohomology classes in the configuration spaces of
    four-manifolds are constructed in an analogous manner.)
    
    Now, fix a curve $\gamma\subset Y$ and consider the one-parameter 
    family of maps 
    $$h_{t}\colon \ModSW(X_{1},\spinc_{1},\CritMan)\longrightarrow S^{1}$$
    indexed by $t\in (0,1]$
    defined by measuring the holonomy of $A$ around the curve 
    $\{1/t\}\times \gamma\subset \XExt_{1}$. 
    Since the configurations in $[A,\Phi]\in \ModSW(X_{1},\spinc_{1},\CritMan)$ 
    converge exponentially to a stationary solution (see~\cite{MMR}), 
    $h_{t}$
    extends continuously to $t=0$. Now, $h_{0}$ represents 
    $(i_{Y}\circ\Restrict{C})^{*}$ of the one-dimensional class 
    $\mu[\gamma]\in H^{1}(\SWConfigIrr(Y,\spinct))$, 
    while $h_{1}$ represents the restriction (to 
    the moduli space) of the
    one-dimensional 
    class $\mu[\gamma_{1}]\in H^{*}(\SWConfigIrr(X_{1},\spinc_{1}))$, where 
    $\mu_{1}=\{1\}\times Y\subset X_{1}$. A similar discusion applies to 
    the two-dimensional class to
    show that $\Restrict{\CritMan}^{*}U=i_{X_{1}}^{*}U$
    (now we use the connection $A$ to 
    identify the fiber ``at infinity'' with the fiber at some point 
    inside $X$). This completes the verification of 
    Inclusion~\eqref{eq:Inclusion}. 
    
    Surjectivity of 
    $$i_{Y}^{*}\colon H^{*}(\SWConfigIrr(Y,\spinct);\Z)\cong \Alg(Y)\rightarrow 
    H^{*}(\CritMan;\Z)$$ 
    follows from 
    classical properties of the cohomology of symmetric products 
    $\Sym^{\mult}(\Sigma)$
    (see~\cite{MacDonald}), according to which the cohomology ring is 
    generated by ``symmetrizations'' of the cohomology of $\Sigma$.
    It is then a straightforward verification 
    (which is spelled out 
    in Proposition~6.10 of~\cite{HigherType})
    using the geometric interpretations of the cohomology 
    classes given  
    above to see that
    that $i_{Y}^{*}\mu(\gamma)$ 
    corresponds to the symmetrization of $\pi(\gamma)$, while 
    $i_{Y}^{*}U$ corresponds to the symmetrization of the point 
    $\pi(x)$ on $\Sigma$ (where we think of $U$ as the Poincar\'e dual 
    of $V_{(x,\Lambda)}$ for some choice of line $\Lambda\subset 
    W_{x}$). 
  
    Thus, it remains to prove that $[\ModSW(X_{1},\spinc_{1},\CritMan)]$ pairs 
    trivially with classes in 
    $$i_{X_{1}}^{*}H^{*}(\SWConfigIrr(X_{1},\spinc_{1})).$$
    We can think of the cohomological pairing $\langle 
    \rho_{_{\CritMan}}^{*}(b), [\ModSW(X_{1},\spinc_{1},\CritMan)]\rangle$ as 
    counting the (signed) number of points to the Seiberg-Witten 
    equations, which satisfy constraints in the compact subset
    $X_{1}\subset\XExt_{1}$; i.e. if $b=U^{d}\cdot 
    [\mu_{1}]\cdot\ldots\cdot[\mu_{\ell}]$, 
    where $\mu_{i}$ are curves in $X_{1}$, and $x_{1},\ldots,x_{d}$ 
    are generic points in $X_{1}$, and 
    $\Lambda_{1},\ldots,\Lambda_{d}$ are generic lines $\Lambda_{i}\subset 
    W^{+}|_{\{x_{i}\}}$, then we have gemetric  
    representatives $V_{(x_{i},\Lambda_{i})}$ and $V_{\mu_{i}}$ for these cohomology 
    classes, so that
    $$\langle [\ModSW(X_{1},\spinc_{1})],b\rangle = \#\ModSW(X_{1})\cap 
    V,$$
    where $V=V_{(x_{1},\Lambda_{1})}\cap\ldots\cap 
    V_{(x_{d},\Lambda_{d})}\cap V_{\mu_{1}}\cap \ldots \cap V_{\mu_{\ell}}$. 
    In fact, if we consider the solutions $\ModSW(X_{1},\spinc_{1},\Jac)$ 
    which satisfy these same constraints, then we get a manifold of 
    dimension one with ends corresponding to $\ModSW(X_{1},\spinc_{1},\CritMan)
    \cap 
    V$. Thus, counting boundary points with sign, we see that
    $$\#\ModSW(X_{1},\spinc_{1},\CritMan)\cap V=0.$$
\end{proof}

\vskip.5cm
\noindent{\bf{Proof of Theorem~\ref{thm:VanishingTheorem}}.} 
In the splitting $X=X_{1}\cup_{Y}X_{2}$, we can number the sides so 
that
the boundary of $X_{1}$ is $Y$ oriented as in Theorem~\ref{thm:ModY}.
Let $\spinc$ be a $\SpinC$ structure over $X$, and let 
$\spinc_{i}\in\SpinC(X_{i})$ denote the restriction 
$\spinc_{i}=\spinc|_{X_{i}}$ for $i=1,2$, and let $\spinct\in\SpinC(Y)$ denote the 
restriction $\spinct=\spinc|_{Y}$.
Let $X(T)$ denote the Riemannian structure on $X$ obtained by inserting
a cylinder $[-T,T]\times Y$ between $X_{1}$ and $X_{2}$ (but keeping 
the metrics on these two pieces to be fixed, and product-like near the 
boundary). If the Seiberg-Witten 
invariants for a $\SpinC$ structure $\spinc$ over $X$ is non-trivial, 
for any unbounded, increasing sequence of real numbers $\{T_{i}\}$, 
there must be a sequence of Seiberg-Witten monopoles 
$[A_{i},\Phi_{i}]\in\ModSW(X(T_{i}),\spinc)$.
The uniform bound in 
the energy, and local compactness (see~\cite{KMthom}) 
allows one to find a sequence $\{t_{i}\}$ with $t_{i}\leq T_{i}$, 
so that, after passing to a subsequence if necessary, 
$[A_{i},\Phi_{i}]|_{\{t_{i}\}\times Y}$ converges to a stationary 
solution; i.e. it converges to a point in $\ModSWThree(Y,\spinct)$. 
Thus, from Theorem~\ref{thm:ModY}, it follows that the Seiberg-Witten 
invariant for $\spinc$ vanishes unless
$\spinct$ is one of the $n$ distinguished $\SpinC$ 
structures $\spinct_{e}$ over $Y$.

Suppose that $\spinct=\spinct_{e}$, for $e=0,\ldots,2g$. 
Note the excision principle for the index gives that
$$d(\spinc)=\dim 
\ModSW(X_{1},\spinc_{1},\CritMan)+\dim\ModSW(X_{2},\spinc_{2},\CritMan) - 
\dim(\CritMan)$$
for any 
$\SpinC$ structure $\spinc\in \SpinC(X)$ with 
$\spinc|_{X_{i}}=\spinc_{i}$ for $i=1,2$ (and generic, compactly 
supported perturbations of the equations over $X_{i}$). We fix an 
integer
$\ell\geq 0$ and homology classes $a_{1},\ldots,a_{m}\in H_{1}(X_{1};\Z)$, 
$b_{1},\ldots,b_{n}\in H_{1}(X_{2};\Z)$ with
$$2\ell+m+n=d(\spinc).$$

Let $\spinc_{1}\#\spinc_{2}\subset \SpinC(X)$ denote the subset of 
$\SpinC$ structures on $X$:
$$\spinc_{1}\#\spinc_{2}=\{\spinc'\in\SpinC(X)\big| 
\spinc'|_{X_{1}}=\spinc_{1}, \spinc'|_{X_{2}}=\spinc_{2}\},$$ and
let
$\ModSW(X(T),\spinc_{1}\#\spinc_{2})$ denote the union
$$\ModSW(X(T),\spinc_{1}\#\spinc_{2})=
\coprod_{\ProdSpinCs}
\ModSW(X,\spinc').$$
Clearly, we have that 
$$\# \ModSW(X(T),\spinc_{1}\#\spinc_{2})\cap V_{1}\cap 
    V_{2} =
    \sum_{\ProdSpinCs}
\SW_{X,\spinc'}(U^{\ell}\cdot[a_{1}]\cdot\ldots\cdot[a_{m}]
\cdot[b_{1}]\cdot\ldots\cdot[b_{n}])
$$
where $V_{i}$ are the intersection of the constraints from the 
$X_{i}$ side; e.g.
$$V_{1}=V_{(x_{1},\Lambda_{1})}\cap\ldots\cap 
V_{(x_{\ell},\Lambda_{m})}
\cap V_{a_{1}}\cap\ldots\cap V_{a_{m}}$$
and 
$$V_{2}=V_{b_{1}}\cap\ldots\cap V_{b_{n}}.$$
Thus, our aim is to prove that the total signed number of points in 
the cut-down moduli space $\ModSW(X(T),\spinc_{1}\#\spinc_{2})\cap 
V_{1}\cap V_{2}$ is zero.

Given pre-compact sets $K_{i}\subset 
\ModSW(X_{i},\spinc_{i},\CritMan)$ for $i=1,2$, there are gluing maps
defined for all sufficiently large $T$, 
$$    \glue_{\CritMan;T}\colon 
K_{1}\#_{\CritMan}K_{2}
\longrightarrow
	\ModSW(X(T),\spinc_{1}\#\spinc_{2}),
$$
where the domain is the fibered product of $K_{1}$ and $K_{2}$
over $\rho_{1}$ and 
$\rho_{2}$, i.e. the set of $[A_{1},\Phi_{1}]\in 
K_{1},[A_{2},\Phi_{2}]\in K_{2}$ with
$$\rho_{1}([A_{1},\Phi_{1}])=
\rho_{2}([A_{2},\Phi_{2}]),$$
and the range consists of all configurations $[A,\Phi]$ which are 
whose restrictions to $X_{1}$ and $X_{2}$ are sufficiently close to 
restrictions (to $X_{1}$ and $X_{2}$)
of configurations $[A_{1},\Phi]\times [A_{2},\Phi_{2}]$ in the 
fibered product. 

We claim that for all sufficiently 
large $T$, the cut-down moduli space 
lies in the range of this map. Specifically, if we had a sequence 
$[A_{i},\Phi_{i}]\in \ModSW(X(T_{i}),\spinc_{1}\#\spinc_{2})$ for an 
increasing, unbounded sequence $\{T_{i}\}_{i=1}^{\infty}$ of real 
numbers, the sequence converges $\Cinfty$ locally to give a pair of 
Seiberg-Witten monopoles
monopoles $[A_{1},\Phi_{1}]\in\ModSW(X_{1},\spinc_{1})$ and 
$[A_{2},\Phi_{2}]\in \ModSW(X_{2},\spinc_{2})$. These monopoles have 
finite energy (since the total variation of $\CSD$ is bounded in the 
limit), so they have boundary values, which must lie in either 
$\CritMan$ or $\Jac$. We exclude all but one of the four cases as 
follows:
\begin{list}
{(C-\arabic{bean})}{\usecounter{bean}\setlength{\rightmargin}{\leftmargin}}
\item The case where $\Restrict{1}[A_{1},\Phi_{1}]\in\Jac$ and 
$\Restrict{2}[A_{2},\Phi_{2}]\in\CritMan$ is excluded since 
$\CSD(\CritMan)<\CSD(\Jac)$. 
\item 
\label{case:RedToRed}
The case where $\Restrict{1}[A_{1},\Phi_{1}]\in\Jac$ and 
$\Restrict{2}[A_{2},\Phi_{2}]\in \Jac$ is excluded by a dimension count. 
Specifically, we must have that
$$\Restrict{1}[A_{1},\Phi_{1}]=\Restrict{2}[A_{2},\Phi_{2}]$$
and $[A_{1},\Phi_{1}]\in \ModSW(X_{1},\Jac)\cap V_{1}$ and 
$[A_{2},\Phi_{2}]\in\ModSW(X_{2},\Jac)\cap V_{2}$, i.e. the pair 
$[A_{1},\Phi_{1}]\times [A_{2},\Phi_{2}]$ lies in the fibered 
product 
$\ModSW(X_{1},\spinc_{1},\Jac)\times_{\Jac}\ModSW(X_{2},\spinc_{2},\Jac)$, 
a space whose dimension is one less than the expected dimension 
$d(\spinc)$ of the 
moduli space. Thus, for generic representatives $V_{1}$ and $V_{2}$, 
this intersection is empty.
\item
\label{case:IrrToRed}
The case where $\Restrict{1}[A_{1},\Phi_{1}]\in\CritMan$ and 
$\Restrict{2}[A_{2},\Phi_{2}]\in\Jac$ is excluded by a similar dimension count. 
We have that $[A_{1},\Phi_{1}]\in\ModSW(X_{1},\spinc_{1},\CritMan)\cap V_{1}$,
$[A_{2},\Phi_{2}]\in\ModSW(X_{2},\spinc_{2},\Jac)\cap V_{2}$,
and $\Restrict{1}[A_{1},\Phi_{1}]$ is connected to $\Restrict{2}[A_{2},\Phi_{2}]$ by 
a (uniquely determined) flow in $\UnparModFlow(\CritMan,\Jac)$. This 
set has expected dimension $-2$. 
\end{list}
The remaining case is that 
$[A_{1},\Phi_{1}]\in\ModSW(X_{1},\spinc_{1},\CritMan)$, and
$[A_{2},\Phi_{2}]\in\ModSW(X_{2},\spinc_{2},\CritMan)$, with
$$\Restrict{1}[A_{1},\Phi_{1}]=\Restrict{2}[A_{2},\Phi_{2}].$$ In particular, 
$[A_{1},\Phi_{1}]$ lies in the compact set 
$\ModSW(X_{1},\spinc_{1},\CritMan)$, while
$[A_{2},\Phi_{2}]$ lies in the set 
$\Restrict{2}^{-1}(\Restrict{1}\ModSW(X_{1},\spinc_{1},\CritMan)\cap V_{1})\cap 
V_{2}$ which is also compact (according to the dimension count used 
to exclude Case~(C-\ref{case:IrrToRed}) above). Thus, for all sufficiently large $T$, the cut-down moduli 
space lies in the image of the gluing map $\glue_{\CritMan;T}$. 

On compact subsets of $X(T)$ away from the ``neck'', 
gluing is a $\Cinfty$ small perturbation, which 
goes to zero as the neck-length is increased; in particular, for 
$i=1,2$,
$$\lim_{T\goesto\infty}\glue_{\CritMan;T}([A_{1},\Phi_{1}]\#[A_{2},\Phi_{2}])
|_{X_{i}}=
[A_{i},\Phi_{i}]|_{X_{i}}.$$
It follows from this that
\begin{equation}
    \label{eq:SWInvariant}
    \#\ModSW(X(T),\spinc_{1}\#\spinc_{2})\cap V_{1}\cap V_{2} =
\#\left(\left(\ModSW(X_{1},\spinc_{1},\CritMan)\cap 
V_{1}\right)\times_{\CritMan}\left(\ModSW(X_{2},\spinc_{2},\CritMan)\cap 
V_{2}\right)\right).
\end{equation}
The latter quantity can be thought of as cohomological pairing in 
$\ModSW(X_{1},\spinc_{1},\CritMan)$ as follows.

Fix an
oriented, $v$-dimensional submanifold $V\subset 
\ModSW(X_{2},\spinc_{2},\CritMan)$, and consider the function which assigns to each 
smooth map
$$f\colon Z\longrightarrow C$$ 
(where $Z$ is a smooth, oriented, compact manifold whose dimensionl 
equals the codimension of $V$)
the number 
of points in the fibered product $\#(Z\times_{\CritMan}V)$ (counting 
with sign, after 
arranging $f$ to be transverse to $V$). 
This is the pairing of the fundamental cycle of $Z$ with an induced 
cohomology class in $H^{d_{2}-v}(\CritMan,\Z)$. Indeed, this class can 
be thought of as the ``push-forward'' of the Poincar\'e dual to $V$, 
under a map
$$(\Restrict{2})_{*}\colon H^{i}(\ModSW(X_{2},\CritMan);\Z)\longrightarrow 
H^{i+\dim(\CritMan)-d_{2}}(\CritMan;\Z).$$
Thus, the count in Equation~\eqref{eq:SWInvariant} can be thought of 
as the pairing
$$\#\ModSW(X(T),\spinc_{1}\#\spinc_{2})\cap V_{1}\cap V_{2} =
\langle [\ModSW(X_{1},\spinc_{1},\CritMan)], \PD(V_{1})\cup 
\Restrict{1}^{*}(\Restrict{2})_{*}\PD(V_{2})\rangle.$$ This pairing vanishes, 
according to Proposition~\ref{prop:TrivialRelativeInvariant}. This 
completes Theorem~\ref{thm:VanishingTheorem} in the case when 
$\spinct=\spinct_{e}$ for $e=0,\ldots,2g-2$.

In the case when $\spinct=\spinct_{e}$ for $2g-1<e<n$,  the 
vanishing of the Seiberg-Witten invariant for any $\spinc$ structure 
with $\spinc|_{Y}=\spinct$ is guaranteed by the same dimension count 
which we used to exclude Case~(C-\ref{case:RedToRed}) above.
\qed

The proof of Theorem~\ref{thm:GenType} follows from an application of
Theorem~\ref{thm:VanishingTheorem}, together with the known properties
of Seiberg-Witten invariants for complex surfaces of general type (see
for instance~\cite{Brussee} or \cite{FriedMorg}, and also~\cite{FS}),
according to which a minimal surface of general type has only two
``basic classes'' ($\SpinC$ structures for which the Seiberg-Witten
invariant is non-zero), the ``canonical'' $\SpinC$ structure
$\spinc_{0}$ (whose first Chern class is given by
$c_{1}(\spinc_{0})=-K_{X}$, where $K_{X}\in H^{2}(X;\Z)$ is the first
Chern class of the complex cotangent bundle of $X$), and its
conjugate. Moreover, the basic classes of the $n$-fold blow-up
${\widehat X}=X\#n\CPbar$ are those $\SpinC$ structures $\spinc$ whose
restriction away from the exceptional spheres agrees with $\spinc_{0}$
or its conjugate, and whose first Chern class evaluated on each
exceptional sphere $E_i$ 
satisfies $$\langle c_{1}(\spinc),[E_{i}]\rangle =
\pm 1.$$ In fact, since $K_{X}^{2}>0$ for a minimal surface of general
type, the basic classes are in one-to-one correspondence with their
first Chern classes.  In view of this fact, throughout the following proof,
we label the basic classes of ${\widehat X}$ by their first Chern 
classes.

\vskip.5cm
\noindent{\bf{Proof of Theorem~\ref{thm:GenType}}.} 
The subgroup $\delta H^{1}(Y;\Z)$ partitions $\SpinC(X)$ into orbits,
and Theorem~\ref{thm:VanishingTheorem} states that if $X$ could be
decomposed, then the sum of invariants under each orbit vanishes. Note
moreover that if $Y$ separates $X$, then the intersection form
restricted to the subgroup $\delta H^{1}(Y;\Z)$ is trivial: this is
true because we can represent cohomology classes $[\omega],[\eta]\in
\delta H^{1}(Y;\R)$ by differential form representatives $\omega$ and
$\eta$, with $\omega|_{X_{1}}\equiv 0$ and $\eta|_{X_{2}}\equiv 0$, so
that the representative for $[\omega]\cup[\eta]$, $\omega\wedge\eta$,
vanishes identically.  It follows from this that in each orbit, there
can exist at most two basic classes, for if we had two basic classes
which had the same coefficient in $K_{X}$, then their difference would
have negative square.  Now, suppose that $K_{X}-E_{1}-\ldots-E_{n}$
had another basic class in its orbit. We know that the other basic
class would be of the form
$$-K_{X}+E_{1}+\ldots+E_{a}-E_{a+1}-\ldots-E_{n},$$ after renumbering
the exceptional curves if necessary.  The difference $\Delta$ is
$2(K_{X}-E_{1}-\ldots-E_{a})$, which must have square zero, which 
forces $a>0$ (recall that 
$K_X^2>0$ for a minimal surface of general type). 
Now, consider the basic class
$K_{X}+E_{1}+\ldots+E_{n}$. It, too, can have at most one other
basic class in its orbit, and the difference has the form
$$\Delta'=2(K_{X}+\epsilon_{1}E_{1}+\ldots\epsilon_{n}E_{n}),$$ where we
know that $\epsilon_{1},\ldots,\epsilon_{n}\geq 0$, in particular
$\Delta-\Delta'$ is a non-zero class, which is easily seen to have
negative square. But this contradicts the fact that $\Delta-\Delta'\in
\delta H^{1}(Y;\Z)$.  Thus, it follows that either the basic class
$K_{X}-E_{1}-\ldots-E_{n}$ or
$K_{X}+E_{1}+\ldots+E_{a}-E_{a+1}-\ldots-E_{n}$ is alone in its
$\delta H^{1}(Y;\Z)$ orbit.  But this contradicts the conclusion of
Theorem~\ref{thm:VanishingTheorem}.
\qed

\subsection{Final Remarks}

It is suggestive to compare the formal framework adopted here with
that of equivariant Morse theory. Specifically, the 
``Chern-Simons-Dirac'' operator on $Y$ in the set-up of 
Theorem~\ref{thm:ModY} has precisely two critical manifolds, a 
manifold of reducibles $\Jac$
(consisting of configurations whose
stabilizer in the gauge group is a circle), and a manifold of 
irreducibles $C$
(consisting of configurations whose stabilizers are trivial). From 
the point of view of equivariant cohomology, then, there should be an 
``equivariant Floer homology'', and an analogue of the Bott spectral 
sequence, whose $E_{2}$ term consists of the homology of the 
irreducible critical point set $H_{*}(\Sym^{e}(\Sigma);\Z)$, and 
the $S^{1}$-equivariant homology of the reducible manifold, which is 
given by
$$H^{*}(\CP{\infty}\times \Jac;\Z)\cong 
\Z[U]\otimes_{\Z}\Wedge^{*}H_{1}(Y;\Z).$$

From this point of view,
Proposition~\ref{prop:TrivialRelativeInvariant}, upon which the 
vanishing theorem rests, can be seen then as the calculation of the 
differential in this spectral sequence. 

The equivariant point of view has been stressed by a number of 
researchers in the field, including (especially in the context of gluing along
rational homology 
three-spheres) \cite{KMlect}, \cite{AusBra}, \cite{MarcolliWang}.

\commentable{
\bibliographystyle{plain}
\bibliography{biblio}
}

\end{document}